\documentclass{article} 

\usepackage{amsmath,amsthm}     
\usepackage{graphicx}     
\usepackage{hyperref} 
\usepackage{url}
\usepackage{amsfonts} 

\usepackage{multicol}

\usepackage{tikz-cd} 
\usepackage{subcaption}

\begin{document}

\title{An uncountable number of proofs of Pythagoras' Theorem}

\author{Gaurav Bhatnagar 
\\ 
Ashoka University,\\
 Sonipat, Haryana 131029, India \\
bhatnagarg@gmail.com \\
\\
Sagar Shrivastava\\  
School of mathematics, \\
Tata Institute of Fundamental Research, \\
Homi Bhabha Road, Mumbai 400005, India\\
sagars@math.tifr.res.in}

\maketitle

\begin{abstract} We give an infinite number of proofs of the Pythagorean theorem. Some can be classified as `self-similar proofs'. 
\end{abstract}

  \usetikzlibrary {calc} 
   \usetikzlibrary {math}
\section{Introduction} 
Virtually anyone who is someone has given a proof of the Pythagorean theorem; we recommend two recent histories given by Maor~\cite{Maor2007} and Veljan~\cite{Veljan2000}, and an interesting set of proofs and commentary by 
Alsina and Nelsen~\cite{AN2015}. The classic compendium and a classification of proofs is by Loomis~\cite{loomis1927}, who says:
\begin{quote}
(there are) many purely geometric demonstrations of this famous theorem accessible to the teacher, as well as an unlimited number of proofs based on the algebraic method of geometric investigation.
\end{quote}
Loomis said this in 1927, so we can say that it has long been suspected that there were an infinite number of proofs by the algebraic method of the Pythagorean theorem. 
However, the proofs of Pythagoras' theorem presented by Loomis and others are still finite in number. The objective of this article is to give an infinite set of proofs of this famous theorem by what Loomis calls the `algebraic method'.

We begin by presenting two proofs. 
One of them is perhaps one of the oldest ever; the other we believe is new. Infinite sequences of proofs are obtained by iteration. Taking the limit, we find proofs that  may be called `self-similar proofs' of Pythagoras' theorem. We call it self-similar, because the figure obtained is self-similar, that is, it contains a scaled copy of itself. This type of proof  has been given by Lengvárszky~\cite{ZL2013}.

The reader may enjoy working out the details of the proofs given in Figure~\ref{fig6}. The picture on the left is 
a proof due to Bh\={a}skara (XIIth century)~\cite[p.~4]{Nelsen1993}, and the one on the right the corresponding self-similar proof, due to Lengvárszky~\cite{ZL2013}.

\begin{figure}[h]
\begin{subfigure}[c]{0.45\textwidth}
\centering 
\begin{tikzpicture}[scale=0.33]
\draw 
   (0,0) node[align=left,   below] {$D$}
-- (13,0) node[align=right,   below] {$C$}
-- (13,13) node[align=right,   above] {$B$}
-- (0,13) node[align=left,   above] {$A$} 
-- (0,0);
\draw (13,0) -- (25/13,60/13);
\draw (13,13) -- (109/13,25/13);
\draw (0,13) -- (144/13,109/13);
\draw (0,0) -- (60/13,144/13);

\draw[below] (5,3) node{$a$};
\draw[below] (1,2) node{$b$};
\draw[below] (2,12) node{$b$};
\draw (6,-0.5) node{$c$};
\end{tikzpicture}
\end{subfigure}
\begin{subfigure}[c]{0.45\textwidth}
\centering
 \begin{tikzpicture}[scale=0.5]

   \coordinate (a) at (0,0);
   \coordinate (b) at (9,0);
   \coordinate (c) at (9,9);
   \coordinate (d) at (0,9);
  
  \draw (a)--(b)--(c)--(d)--cycle;

  \foreach \i in {1,2,...,15}
  {

   \coordinate (a1) at ($ (b)!0.95!-18:(a) $);
   \coordinate (b1) at ($ (c)!0.95!-18:(b) $);
   \coordinate (c1) at ($ (d)!0.95!-18:(c) $);
   \coordinate (d1) at ($ (a)!0.95!-18:(d) $);

 \draw (a1)--(b);
 \draw (b1)--(c);
 \draw (d)--(c1);
 \draw (a)--(d1);
   
   \coordinate (a) at (a1);
   \coordinate (b) at (b1);
   \coordinate (c) at (c1);
   \coordinate (d) at (d1);
   };

\end{tikzpicture}

\end{subfigure} 

\caption{Bh\={a}skara's proof of the Pythagorean theorem and the corresponding self-similar proof}\label{fig6}
\end{figure}

In this paper, we motivate and provide two other proofs of this nature. 
Next, we show how to modify our approach to obtain an uncountable number of algebraic proofs. Previously, Alsina and Nelsen~\cite[p.~47]{AN2006} (see also
\cite[p.~257]{AN2015}) have also indicated how to obtain an uncountable number of proofs of Pythagoras' theorem via tiling. Our  approach yields an infinite number of self-similar proofs too---however, these are only countably many in number. 

Finally, we note that we can reverse-engineer one of the proofs, to obtain the formula for the geometric sum---assuming the Pythagorean 
theorem---which leads to a nice surprise for a calculus student.

\section{Two proofs of the Pythagorean theorem}\label{sec2}
We begin with two related proofs of the Pythagorean theorem. The first proof is based on a figure (see Figure~\ref{fig1}) from the earliest known proof of the theorem~\cite[p.~3]{Nelsen1993}. This has been credited to Chou-pei Suan-ching (circa 250 B.C.) by
\cite{Veljan2000}. The second appears to be new. 

\subsection*{Proof \# 1}

We begin with the picture in Figure~\ref{fig1}. 
Consider any right-angled triangle with perpendicular $a_1$ and base $b_1$ with the hypotenuse of length $c_1$. We begin with a square $ABCD$ of side length $a_1+b_1$. Let $E$ be the point on the edge $AB$ such that $|AE|=a_1$ and $|EB|=b_1$. Similarly $F,G,H$ be points on $BC,CD,DA$ 
respectively, with lengths as shown.

\begin{figure}[h]
\centering
\begin{tikzpicture}[scale=0.75]
\draw 
   (0,0) node[align=left,   below] {$D$}
-- (9,0) node[align=right,   below] {$C$}
-- (9,9) node[align=right,   above] {$B$}
-- (0,9) node[align=left,   above] {$A$} 
-- (0,0);
\draw 
   (0,3) -- (6,0) -- (9,6) -- (3,9) -- (0,3);

\draw (-0.3,3) node{$H$};
\draw (6,-0.3) node{$G$};
\draw (9.3,6) node{$ F$};
\draw (3,9.3) node{$ E$};
\draw (1.5,9.25) node{$a_1$};
\draw (6,9.25) node{$b_1$};
\draw (9.3,7.5) node{$a_1$};
\draw (9.3,3) node{$b_1$};
\draw (3,-0.2) node{$b_1$};
\draw (7.5,-0.2) node{$a_1$};
\draw (-0.2,1.5) node{$a_1$};
\draw (-0.2,6) node{$b_1$};

\draw (6,7.8) node{$c_1$};
\draw (1.5,2) node{$c_1$};
\draw (1.5,6.5) node{$c_1$};
\draw (6.7,2) node{$c_1$};
\end{tikzpicture}

\caption{Proof \# 1: Chou-pei Suan-ching's proof.} \label{fig1}
\end{figure}

It is easy to verify (using the fact that the sum of angles in a triangle is $180^\circ$) that the
quadrilateral $EFGH$ is a square.  Now computing the area of the square of side $a_1+b_1$ in two different ways, we have
$$(a_1+b_1)^2 = 4\times \frac{1}{2} a_1b_1+ c_1^2,$$
from where we obtain
$$a_1^2+b_1^2=c_1^2.$$
This completes a proof of the Pythagorean theorem.

\subsection*{Proof \# 2}
Next we modify the picture of the first proof as shown in Figure~\ref{fig2a} to obtain another proof of the Pythagorean theorem. This proof appears to be new.

\begin{figure}[h]
\centering
\begin{tikzpicture}[scale=0.75]
\draw 
   (0,0) node[align=left,   below] {$D$}
-- (9,0) node[align=right,   below] {$C$}
-- (9,9) node[align=right,   above] {$B$}
-- (0,9) node[align=left,   above] {$A$} 
-- (0,0);
\draw 
   (0,3) -- (6,0) -- (9,6) -- (3,9) -- (0,3);

\draw (-0.3,3) node{$H$};
\draw (6,-0.3) node{$G$};
\draw (9.3,6) node{$F$};
\draw (3,9.3) node{$E$};
\draw (2,2) -- (7,2) -- (7,7) -- (2,7) -- (2,2);
\draw (1.9,1.8) node{$K$};
\draw (7.2,1.9) node{$J$};
\draw (7.2,7.2) node{$I$};
\draw (1.8,7.1) node{$L$};
\draw (9,7.5) node[right]{$a_1$};
\draw (9,3) node[right]{$b_1$};
\draw (3,9) --(3,7) node[below] {$M$};
\draw (6,0) --(6,2) node[above] {$N$};
\draw (3,8) node[right]{$d_1$};
\draw (6,1) node[left]{$d_1$};
\draw (8,6.5) node[above]{$a_2$};
\draw (7,4.5) node[left]{$c_2$};
\draw (5,8) node[above]{$b_2$};

\draw (8,4) node[right]{$b_2$};

\end{tikzpicture}

\caption{Proof \# 2.}\label{fig2a}
\end{figure}
Let
$I$ be the point on the edge $EF$ such that $\frac{|EI|}{|IF|}=\frac{b_1}{a_1}$. 
We denote by $a_2$ and $b_2$ the lengths $|IF|$ and $|EI|$. 

Similarly we have points $J,K,L$ on $FG,GH,HE$, such that 
$$\frac{|FJ|}{|JG|}=\frac{|GK|}{|KH|}=\frac{|HL|}{|LE|}=\frac{b_1}{a_1}.$$ 

It is easy to see that the triangles $\bigtriangleup ELI, \bigtriangleup FIJ, \bigtriangleup GJK, \bigtriangleup HKL$ are all right-angled triangles, and congruent to each other. Let
$c_2$ denote the length of the hypotenuse of these triangles. 

The following are easy to prove. 
\begin{itemize}
\item The triangle $\bigtriangleup ELI$  is similar to the bigger triangles from the first proof, namely $\bigtriangleup AEH,$ $ \bigtriangleup BFE, \bigtriangleup CGF, \bigtriangleup DHG$. So too are $\bigtriangleup FIJ$, $\bigtriangleup GJK$ and $\bigtriangleup HKL$.
\item The quadrilateral $IJKL$ is a square. 
\item The edges $AB$ and $LI$ are parallel. 
\item The lengths $a_2$, $b_2$ and $c_2$ are given by
$$
a_2= \frac{a_1c_1}{a_1+b_1} ; \;\;
b_2= \frac{b_1c_1}{a_1+b_1}; \;\;
c_2= \frac{c_1^2}{a_1+b_1} .
$$
\end{itemize}
Let $EM$ be the perpendicular dropped from $E$ on $LI$; similarly $GN$ is the perpendicular dropped from $G$ on $JK$. As these are both perpendiculars dropped in the same orientation of congruent triangles, they have the same length, say $d_1= |EM|=|GN|$. Thus
\begin{equation} \label{eq:1}
a_1+b_1 =|BC| = |EM|+|IJ|+|GN| = 2d_1 +c_2
\end{equation}
To find the value of $d_1$, compute the area of $\bigtriangleup ELI$ in two ways:
\begin{align*}
\operatorname{area}(\bigtriangleup ELI) &= \frac{1}{2} |EL| |EI| = \frac{c_1^2 a_1 b_1}{2(a_1+b_1)^2};\\
\intertext{and,}
\operatorname{area}(\bigtriangleup ELI) &= \frac{1}{2} |LI| |EM| = \frac{c_1^2 d_1}{2(a_1+b_1)}.
\end{align*}
Equating the two gives us $$d_1 = \frac{a_1b_1}{a_1+b_1}.$$

Putting these values in \eqref{eq:1}, we get:
$$a_1+b_1 = \frac{c_1^2}{a_1+b_1} + \frac{2a_1b_1}{a_1+b_1}.$$
Multiplying the equation by $a_1+b_1$ gives us $a_1^2 +b_1^2 =c_1^2$. This completes the second proof of Pythagoras' theorem.

\section{Self-similar proofs of the Pythagorean theorem}

An iterate of the first proof is hidden in Figure~\ref{fig2a}. This leads to an infinite sequence of proofs.

\begin{figure}[h]
\begin{center}
\begin{tikzpicture}[scale=0.75]

\draw 
   (0,0) 
-- (9,0) 
-- (9,9) 
-- (0,9)  
-- (0,0);
\draw 
   (0,3) -- (6,0) -- (9,6) -- (3,9) -- (0,3);

\draw (2,2) -- (7,2) -- (7,7) -- (2,7) -- (2,2);
\draw (1.5,9.3) node{$a_1$};

\draw (7,9) node[above]{$b_1$};

\draw (9,7.5) node[right]{$a_1$};

\draw (2.5,8) node[left]{$a_2$};
\draw (4.5,7) node[above]{$c_2$};
\draw (5,8) node[above]{$b_2$};

\draw (8,6.5) node[above]{$a_2$};
\end{tikzpicture}

\end{center}

\caption{The third proof}\label{fig3}
\end{figure}

Consider Figure~\ref{fig3}.
On the one hand, the area of the larger square is $(a_1+b_1)^2$. On the other hand, from
Figure~\ref{fig3}, we see  that it equals
\begin{align*}
4\times & \frac{1}{2}a_1b_1 + 4\times \frac{1}{2}a_2b_2 + c^2_2
\\
&
=
2a_1b_1+ c_1^2.
\end{align*} 
The last step follows by noting that $c_2^2+2a_2b_2=c_1^2$, a fact that is obvious from the figure.
Now this reduces to the earlier proof indicated by Figure~\ref{fig1}. 

\subsection*{Proof \# n+1}
Now we iterate. Suppose that we make $n$ squares (inside the square of size $a_1+b_1$) by repeating the above steps. Let $a_n$, $b_n$ and $c_n$ be the relevant lengths. 
We find that 
\begin{align}\label{recursive-proof}
(a_1+b_1)^2 & = 2a_1b_1+2a_2b_2+\cdots+2a_{n-1}b_{n-1}+2a_nb_n+ c_n^2 \\
& = 2a_1b_1+2a_2b_2+\cdots+2a_{n-1}b_{n-1}+ c_{n-1}^2 \cr
& \dots\cr
&= 2a_1b_1+c_1^2.\nonumber
\end{align}

In fact, the above gives an infinite (but still countable) set of proofs of Pythagoras' theorem.

However, if we iterate endlessly, we get another proof---this time obtaining the self-similar image given in Figure~\ref{fig4}.

\begin{figure}[h]

\centering
 \begin{tikzpicture}

   \coordinate (a) at (0,0);
   \coordinate (b) at (9,0);
   \coordinate (c) at (9,9);
   \coordinate (d) at (0,9);
  
  \draw (a)--(b)--(c)--(d)--cycle;

  \foreach \i in {1,2,...,15}
  {

   \coordinate (a1) at ($ (a)!0.7!(b) $);
   \coordinate (b1) at ($ (b)!0.7!(c) $);
   \coordinate (c1) at ($ (c)!0.7!(d) $);
   \coordinate (d1) at ($ (d)!0.7!(a) $);

 \draw (a1)--(b1)--(c1)--(d1)--cycle;
   
   \coordinate (a) at (a1);
   \coordinate (b) at (b1);
   \coordinate (c) at (c1);
   \coordinate (d) at (d1);
   };

\end{tikzpicture}
\caption{Self-similar proof \#1 of the Pythagorean theorem}\label{fig4}
\end{figure}

\subsection*{Self-similar proof \#1}
We now begin the proof from Figure~\ref{fig4}. We note that
the lengths $a_n$, $b_n$ and $c_n$ (for $n> 1$) are given by
$$
a_{n+1}= \frac{a_nc_n}{a_n+b_n} ; \;\;
b_{n+1}= \frac{b_nc_n}{a_n+b_n}; \;\;
c_{n+1}= \frac{c_n^2}{a_n+b_n} .
$$
In addition, recall that $c_n=a_{n+1}+b_{n+1}$.
From these, we obtain formulas for $a_n$, $b_n$ and $c_n$, by iteration. Note that 
$$
\frac{a_{n+1}}{a_n}= \frac{b_{n+1}}{b_n}=  \frac{c_{n+1}}{c_n}=  \frac{c_{n}}{a_{n}+b_{n}},
$$
so
$$\frac{c_{n+1}}{c_n}=  \frac{c_{n}}{a_{n}+b_{n}}=\frac{c_n}{c_{n-1}}=\cdots =\frac{c_2}{c_1}= \frac{c_{1}}{a_{1}+b_{1}}.$$
From these we obtain:
\begin{equation}\label{an-bn-cn}
a_n= a_1\Big( \frac{c_1}{a_1+b_1}\Big)^{n-1}; 
b_n= b_1\Big( \frac{c_1}{a_1+b_1}\Big)^{n-1};
c_n= c_1\Big( \frac{c_1}{a_1+b_1}\Big)^{n-1}.
\end{equation}

We know the area of the largest square is $(a_1+b_1)^2$. The other way is to consider the sum of the areas of the individual triangles making up the square. The area of the outermost four triangles is $4  \times {a_1b_1}/{2}= 2a_1b_1$; the area of the next four $2a_2b_2$, and so on. 
The total area is evidently:
\begin{align*}
\text{Area of square} &= 2a_1b_1 + 2a_2 b_2 +2a_3b_3+\cdots \\
 &=2a_1b_1 + 2a_1b_1 \left( \frac{c_1}{a_1+b_1} \right)^2 + 2a_1b_1 \left( \frac{c_1}{a_1+b_1} \right)^4 + \cdots \\
&= \frac{2a_1b_1}{1-(\tfrac{c_1}{a_1+b_1})^2},
\end{align*}
where we have used \eqref{an-bn-cn} and summed the geometric series.
Equating the two areas, we get:
$$(a_1+b_1)^2 = \frac{2a_1b_1}{1-(\tfrac{c_1}{a_1+b_1})^2} \implies (a_1+b_1)^2 - c_1^2 = 2a_1b_1 \implies a_1^2 +b_1^2=c_1^2.$$

This completes the proof. Interestingly, to use the sum of the geometric series we need to know that the sum of any two side lengths of a triangle is bigger than the third, and so $c_1/(a_1+b_1)<1$. 

\subsection*{Self-similar proof \#2}
If we iterate Figure~\ref{fig2a}, we will get another countable sequence of proofs, and in the limit, a self-similar proof of Pythagoras' theorem (see Figure~\ref{fig5}). It is easy to see that the iterates of $d_1$ satisfy
$$d_{n} = \frac{a_{2n-1}b_{2n-1}}{a_{2n-1}+b_{2n-1}}.$$
The proof is obtained from
$$a_1+b_1=2d_1+2d_2+2d_3+\cdots.$$

\begin{figure}[h]
\centering
 \begin{tikzpicture}

   \coordinate (a) at (0,0);
   \coordinate (b) at (9,0);
   \coordinate (c) at (9,9);
   \coordinate (d) at (0,9);
   
  \draw (a)--(b)--(c)--(d)--cycle;

  \foreach \j in {0,1}
  {

              \coordinate (a1) at ($ (a)!0.7!(b) $);
            \coordinate (b1) at ($ (b)!0.7!(c) $);
            \coordinate (c1) at ($ (c)!0.7!(d) $);
            \coordinate (d1) at ($ (d)!0.7!(a) $);

      \draw (a1)--(b1)--(c1)--(d1)--cycle;
   
      \coordinate (a) at (a1);
      \coordinate (b) at (b1);
      \coordinate (c) at (c1);
      \coordinate (d) at (d1);

             \coordinate (a1) at ($ (a)!0.3!(b) $);
            \coordinate (b1) at ($ (b)!0.3!(c) $);
            \coordinate (c1) at ($ (c)!0.3!(d) $);
            \coordinate (d1) at ($ (d)!0.3!(a) $);

      \draw (a1) -| (a);
      \draw (c1) -| (c);
      
      \draw (a1)--(b1)--(c1)--(d1)--cycle;
   
      \coordinate (a) at (a1);
      \coordinate (b) at (b1);
      \coordinate (c) at (c1);
      \coordinate (d) at (d1);

                  \coordinate (a1) at ($ (a)!0.7!(b) $);
            \coordinate (b1) at ($ (b)!0.7!(c) $);
            \coordinate (c1) at ($ (c)!0.7!(d) $);
            \coordinate (d1) at ($ (d)!0.7!(a) $);

      \draw (a1)--(b1)--(c1)--(d1)--cycle;
   
      \coordinate (a) at (a1);
      \coordinate (b) at (b1);
      \coordinate (c) at (c1);
      \coordinate (d) at (d1);

             \coordinate (a1) at ($ (a)!0.3!(b) $);
            \coordinate (b1) at ($ (b)!0.3!(c) $);
            \coordinate (c1) at ($ (c)!0.3!(d) $);
            \coordinate (d1) at ($ (d)!0.3!(a) $);

      \draw (b1) -| (b);
      \draw (d1) -| (d);
      
      \draw (a1)--(b1)--(c1)--(d1)--cycle;
   
      \coordinate (a) at (a1);
      \coordinate (b) at (b1);
      \coordinate (c) at (c1);
      \coordinate (d) at (d1);
  };
  
\end{tikzpicture}

\caption{Self-similar proof \#2 of the Pythagorean theorem}\label{fig5}
\end{figure}

\subsection*{Self-similar proof \#3}
Another classical proof of Pythagoras' theorem uses the figure on the left in Figure~\ref{fig6} in the introduction. We leave it to the reader to complete the details of the self-similar proof (given on the right) arising from this proof. The details are similar to those in the proofs above.

\section{Uncountably many proofs of the Pythagorean theorem}
There are uncountably many pictures which can be obtained from variations of 
Figure~\ref{fig4} which all lead to a proof of the Pythagoras' theorem. At each step, one can place the inner square in two ways. We associate a binary string with each figure. The digits $0$ and $1$ are attached to the figures shown in Figure~\ref{fig7}.

\begin{figure}[h]

\begin{center}
\begin{tikzpicture}[scale=0.75]
\draw 
   (0,0) 
-- (5,0) 
-- (5,5) 
-- (0,5)  
-- (0,0);
\draw 
   (2,5) -- (5,3) -- (3,0) -- (0,2) -- (2,5);

\draw 
   (7,0) 
-- (12,0) 
-- (12,5) 
-- (7, 5)  
-- (7,0);
\draw 
   (10,5) -- (12,2) -- (9,0) -- (7,3) -- (10,5);

\draw (1,5) node[above]{$a_n$};
\draw (3.5,5) node[above]{$b_n$};

\draw (8.5,5) node[above]{$b_n$};
\draw (10.5,5) node[above]{$a_n$};

\draw (2.5,0) node[below]{$\begin{matrix}\updownarrow\\ 1\end{matrix}$};
\draw (9.5,0) node[below]{$\begin{matrix}\updownarrow\\ 0\end{matrix}$};

\end{tikzpicture}

\end{center}

\caption{The correspondence of binary digit with the orientation of square}\label{fig7}
\end{figure}

Thus Figure~\ref{fig1} is associated to the binary string $1$, Figure~\ref{fig3} to $10$ and Figure~\ref{fig4} to $000000\dots$.  Figure~\ref{fig8} shows the figure corresponding to the binary string
$010101$.

\begin{figure}[h]
\centering
 \begin{tikzpicture}[scale=0.75]

   \coordinate (a) at (0,0);
   \coordinate (b) at (9,0);
   \coordinate (c) at (9,9);
   \coordinate (d) at (0,9);
  \draw (a)--(b)--(c)--(d)--cycle;

\foreach \j in {0,1,0,1,0,1}
  {
  \tikzmath{
  real \x;
  \i=0.3+0.4*\j;
  }
            \coordinate (a1) at ($ (a)!\i!(b) $);
            \coordinate (b1) at ($ (b)!\i!(c) $);
            \coordinate (c1) at ($ (c)!\i!(d) $);
            \coordinate (d1) at ($ (d)!\i!(a) $);

      \draw (a1)--(b1)--(c1)--(d1)--cycle;
   
      \coordinate (a) at (a1);
      \coordinate (b) at (b1);
      \coordinate (c) at (c1);
      \coordinate (d) at (d1);
   };
  
\end{tikzpicture}

\caption{The figure corresponding to $010101$}\label{fig8}

\end{figure}

Notice that in the correspondence with binary numbers, whenever we have a binary number which repeats a sequence of digits, the corresponding figure is a  self-similar diagram. This kind of binary expansion corresponds to a rational number. In general, the binary expansion corresponding to a rational number either terminates (and has an infinite string of $0$s at the end) or some sequence of digits eventually repeat. In either case, the corresponding diagram contains a self-similar figure.  On the other hand, if it's a diagram corresponding to the binary expansion of an irrational number, then there can't be any repeating patterns, as that would contradict the irrational nature of the number. Thus there are only countably many self-similar proofs obtainable by this method!

We remark that there can be further uncountable families of such proofs, two being provided by 
analogous variations in Figures~\ref{fig6} and \ref{fig5}.

\section{Pythagoras' theorem implies the geometric sum}
This paper has many proofs of the Pythagorean theorem. If we assume the Pythagorean theorem, we can prove the formula for the sum of the geometric series, by reverse-engineering Proof \# $n+1$. 

From equations 
\eqref{recursive-proof} and \eqref{an-bn-cn} we have
\begin{multline*}
(a_1+b_1)^2 = 2a_1b_1+2a_1b_1\Big( \frac{c_1}{a_1+b_1}\Big)^{2}+\cdots \\ + 
2a_1b_1\Big( \frac{c_1}{a_1+b_1}\Big)^{2n-2}+ c_1^2 \Big( \frac{c_1}{a_1+b_1}\Big)^{2n-2} 
\end{multline*}
or 
$$\frac{(a_1+b_1)^2 - {c_1^{2n}}/{(a_1+b_1)^{2n-2}  }} {2a_1b_1} = 
1+\Big( \frac{c_1}{a_1+b_1}\Big)^{2}+\cdots+\Big( \frac{c_1}{a_1+b_1}\Big)^{2n-2}
.$$
By Pythagoras' theorem, $ 2a_1b_1=(a_1+b_1)^2 - c_1^2 $, and this expression can be written as
$$
\frac{1 - \big(\frac{c_1}{a_1+b_1 }\big)^{2n}}  {1-\big(\frac{c_1}{a_1+b_1 }\big)^{2} } = 
1+\Big( \frac{c_1}{a_1+b_1}\Big)^{2}+\cdots+\Big( \frac{c_1}{a_1+b_1}\Big)^{2n-2}
.$$

We have obtained the formula for the geometric sum, but with the restriction that $a_1$, $b_1$ and $c_1$ form the sides of a right-angled triangle. That restriction can be removed by considering it once again as an algebraic identity. Replace $c_1$ by $\sqrt{q} (a_1+b_1)$, to obtain the geometric sum:
\begin{equation}\label{geometric}
1+q+q^2+\cdots + q^{n-1} = \frac{1-q^n}{1-q}.
\end{equation}
We can argue that \eqref{geometric} holds for all real values of $q$ provided $q\neq 1$. The argument is as follows. Multiply both sides of \eqref{geometric} by $(1-q)$, to find that both sides are polynomials in $q$ of degree $n$. The polynomials are equal for infinitely many values of $q$ (for $1/2 \le q < 1$), and thus the equation is an identity. On dividing by $(1-q)$, we obtain the formula for the geometric sum. 

Calculus students will recognize in this a way to find the left derivative of $f(q)=q^n$ at $q=1$, by taking the limit as $q\to 1^-$; or, in other words from the 
Pythagorean theorem by taking $b_1\to 0$, so that
$a_1\to c_1$ and thus $q\to 1^-$. Geometrically, this corresponds to squashing the right-angled triangle until it becomes a straight line of length $c_1$.

\subsection*{Acknowledgements} 
We thank Roger Nelsen and Tom Edgar for bringing some references to our attention. We are grateful to the referee for greatly improving the exposition.

\scriptsize

\noindent
{\bfseries GAURAV BHATNAGAR} grew up in Delhi and obtained his Ph.D. at the Ohio State University.  He is most interested in Ramanujan's identities and spends much of his time trying to understand and explain how Ramanujan could possibly have discovered them. In his research, he has discovered and proved several countably infinite families of identities.
\vspace {10pt}

\noindent
{\bfseries SAGAR SHRIVASTAVA} is completing his Ph.D. at the School of Mathematics, Tata Institute of Fundamental Research (TIFR), Mumbai.
His research is in the area of Representation Theory. He participates in mathematics outreach programs where he uses Origami to visualize hyperbolic and planar geometrical surfaces, and to construct platonic solids. During his studies, he has encountered a seemingly infinite number of theorems and proofs.


\begin{thebibliography}{1}
	
	\bibitem{AN2006}
	C.~Alsina and R.~B. Nelsen.
	\newblock {\em Math made visual: {C}reating images for understanding
		mathematics}, volume~28.
	\newblock American Mathematical Soc., 2006.
	
	\bibitem{AN2015}
	C.~Alsina and R.~B. Nelsen.
	\newblock {\em Icons of mathematics. {An} exploration of twenty key images},
	volume~45 of {\em Dolciani Math. Expo.}
	\newblock Washington, DC: Mathematical Association of America (MAA), 2011.
	
	\bibitem{ZL2013}
	Z.~Lengvárszky.
	\newblock Proving the {P}ythagorean theorem via infinite dissections.
	\newblock {\em The American Mathematical Monthly}, 120(8):751--753, 2013.
	
	\bibitem{loomis1927}
	E.~S. Loomis.
	\newblock The {Pythagorean} proposition, its proofs analyzed and classified,
	and bibliography of sources.
	\newblock 214 p. {Ohio}, 1927.
	
	\bibitem{Maor2007}
	E.~Maor.
	\newblock {\em The {Pythagorean} theorem. {A} 4000-year history}.
	\newblock Princeton Sci. Libr. Princeton, NJ: Princeton University Press,
	reprint of the 2007 hardback edition edition, 2019.
	
	\bibitem{Nelsen1993}
	R.~B. Nelsen.
	\newblock {\em Proofs without words. {Exercises} in visual thinking}, volume~1
	of {\em Classr. Resour. Mater.}
	\newblock Washington, DC: Mathematical Association of America (MAA), 1993.
	
	\bibitem{Veljan2000}
	D.~Veljan.
	\newblock The 2500-year-old {P}ythagorean theorem.
	\newblock {\em Mathematics Magazine}, 73(4):259--272, 2000.
	
\end{thebibliography}
\end{document}